\documentclass{ws-apjor}
\usepackage{amsmath,amsfonts,amssymb,verbatim}
\usepackage[round]{natbib}

 \newtheorem{Rem}{Remark}

 \newtheorem{Def}{Definition}
 \newtheorem{Qu}{Question}

 \def\beQ{\begin{Qu}}  \def\eeQ{\end{Qu}}
 \def\beR{\begin{Rem}} \def\eeR{\end{Rem}}

\newcommand{\vc}[1]{{\boldsymbol #1}}

 \def\bep{\begin{pmatrix}}  \def\eep{\end{pmatrix}}

    \def\BI{\begin{itemize}} \def\EI{\end{itemize}}
    \def\im{\item}

    \def\ssec{\subsection}   \def\sec{\section}
    
     \def\no{\noindent}

     \def\fr{\frac}
    \def\beq{\begin{eqnarray}} \def\eeq{\end{eqnarray}}

 \def\bea{\begin{eqnarray*}} \def\eea{\end{eqnarray*}}
 \def\be{\begin{enumerate}} \def\ee{\end{enumerate}}
 \def\im{\item}

 \def\I{\infty}

 \def\a{\alpha } 
\def\b{\beta}    \def\vp{\varphi}

 \def\l{\lambda}

\def\th{\theta}  \def\Tth{\T{\theta}} 
 \def\le{\left} \def\ri{\right}

\def\1{\mathop{\rm 1\!\!I}\nolimits}

\def\C{{\mathbb C}} 
 
 \def\N{\mathbb N}
 \def\le{\left} \def\ri{\right} \def\I{\infty}
   
\def\N{{\mathbb N}}

  \def\T{\tilde}

\def\q{\quad} \def\fr{\frac}
\def\Bmu{\bar \mu}

\def\la{\label}

 \def\no{\nonumber}

    \newcommand{\bff}[1]{{\mbox{\boldmath$#1$}}}
\def\Bmu{\bar \mu}

\def\sec{\section} \def\ssec{\subsection}

\long\def\symbolfootnote[#1]#2{
\begingroup
\def\thefootnote{\fnsymbol{footnote}}\footnote[#1]{#2}
\endgroup}
\def\fn{\symbolfootnote}
\def\Tl{\tilde{\lambda}}

   \def\Tl{\tilde{\lambda}}
\def\Tr{\tilde{\rho}}

\def\Tm{\tilde{\m}} \def\Tmu{\tilde{\mu}}
  
\def\x{\xi}    \def\r{\rho} \def\s{\sigma}
 \def\k{\kappa}
\def\c{s}
\def\n{\nu} \def\m{\mu}  %\def\n{\mu_r}

\def\la{\label}

 \def\l{\lambda} \def\r{\rho}

\def\th{\theta}  
\def\ri{\right}

\def\A{A%^{(0)}
} \def\B{B%^{(1)}
} \def\C{C%^{(2)}
}
\def\j{n}

\begin{document}

%\markboth{Authors' Names}{Instructions for
%Typing Manuscripts (Paper's Title)}

%%%%%%%%%%%%%%%%%%%%% Publisher's Area please ignore %%%%%%%%%%%%%%%
%
%\catchline{}{}{}{}{}
%
%%%%%%%%%%%%%%%%%%%%%%%%%%%%%%%%%%%%%%%%%%%%%%%%%%%%%%%%%%%%%%%%%%%%

\title{\uppercase{\bf{On
multiserver
retrial queues:
history, Okubo-type  hypergeometric systems and matrix continued-fractions}}
}

\author{F. AVRAM}

\address{D\'epartement de Math\'ematiques,
Avenue de l'Universit\'e \\Pau, France\\
 florin.avram@univ-pau.fr\\
http://web.univ-pau.fr/~avram/}

\author{D. MATEI}

\address{Institute of Mathematics ``Simion Stoilow'' of the Romanian Academy\\
Bucharest, P.O. Box 1-764, RO-014700, Romania\\
dmatei@imar.ro
}

\author{Y. Q. ZHAO}

\address{School of Mathematics and Statistics,
Carleton University \\
Ottawa, Canada \\
zhao@math.carleton.ca \\
http://people.math.carleton.ca/~zhao/}

\maketitle

\begin{history}
 \received{(17/07/2012)}
 \revised{(4/08/2013)}
%\accepted{(Day Month Year)}
%\comby{(xxxxxxxxxx)}
\end{history}

\begin{abstract}
%The abstract should summarize the context, content
%and conclusions of the paper in less than 200 words. It should
%not contain any reference citations or displayed equations. Typeset
%the abstract in 8 pt roman with baselineskip of 10 pt, making
%an indentation of 1.5 pica on the left and right margins.

We study two families of QBD processes with linear rates: (A) the
multiserver retrial queue and its easier relative; and (B) the
multiserver M/M/$\I$ Markov modulated queue.

The linear rates imply that  the stationary probabilities satisfy
a recurrence with linear coefficients;  as known from previous
work, they yield a  ``minimal/non-dominant" solution of this
recurrence, which may be computed numerically by matrix
continued-fraction methods.

Furthermore, the generating function of the stationary
probabilities satisfies a linear differential system with
polynomial coefficients, which calls for the venerable but still
developing theory of holonomic (or D-finite) linear differential
systems.

We provide  a  differential system for our  generating function
that unifies problems (A) and (B), and we also include some
additional features and   observe  that in at least one particular
case we get  a special ``Okubo-type hypergeometric system",  a
family that recently spurred considerable interest.

The differential system should allow  further study of the Taylor
coefficients of the expansion of the generating function at  three
points of interest: 1) the irregular singularity at $0$; 2) the
dominant regular singularity, which yields  asymptotic series via
classic methods like the Frobenius vector expansion;  and 3) the
point $1$, whose Taylor series coefficients are the factorial moments.
\end{abstract}


\sec{Introduction} \vskip.1truecm{\bf Motivation:} The practical
questions of taking into account the blocking of  queues, the
dilemma of blocked customers either waiting or  retrying later,
and the eventual abandon  of
 dissatisfied customers have been recognized by queueing theorists from the outstart
-- see for example \citet{kosten1947influence}.

The modelling of these phenomena is  more challenging
mathematically than that of classical  queues, since even the
simplest M/M/1 + retrial case involves a bivariate modelling that
tries to capture the interaction between two types of customers:
those  {who ``}wait on the spot", and those who
 {``}postpone their  {requests} for later".

\vskip.1truecm{\bf A bit of history: the mystery of retrial queues
with more than two servers.} An  impressive pioneering analysis of
retrial queues stems from J.W. Cohen \citeyearpar{cohen1957basic},
who obtained  exact solutions in terms of contour integrals and
Laguerre series expansions. For $\c = 1$ and $\c= 2$ servers,
Cohen's results involve classic hypergeometric functions.

One of the striking features of the field is the {``}disappearing"
of second-order classic hypergeometric equations for the
generating function of the stationary probabilities  when the
number $s$ of servers is bigger than $2$. This has been noted
several times in the  literature -- see, for example pp.25 of
\citet{kulkarni1997retrial} and pp.288 of
\citet{falin1997retrial}. Also, in the case of  $\c = 3$ servers,
there is a considerable increase in complexity that forces one to
turn to numerical methods  as  was  done by \citet{kim95retrial}
and G\'{o}mez-Corral and
Ramalhoto~\citeyearpar{gomez1999stationary}.

\vskip.1truecm{\bf The non-Fuchsian singularity}. Let us dwell
further on this important type of singularity.   The generating
function of the stationary probability distribution satisfies the
well-known Kolmogorov linear differential system. In our case,
this has affine  coefficients --- see for example
\citet{riordan1962stochastic}, \citet{keilson1968service},
\citet{cohen1957basic}, \citet{falin1997retrial},
\citet[(2.2)]{artalejo2002standard},  and our generalization
\eqref{linsys}.

With a low number of servers ($\c \leq 2$), the system leads to a
regular singularity at $0$  (which explains the explicit
hypergeometric solutions), while with a higher number of servers,
it has a non-regular (non-Fuchsian) singularity.

\vskip.1truecm{\bf Formal series solutions}. The concept of
classic (second order) hypergeometric functions has been
fruitfully extended to that of functions satisfying holonomic
(D-finite) equations, for example see
\citet{stanley1980differentiably,zeilberger1990three,koepf1997algebra}.
The  {\bf general solution}  of such equations
 falls straight under the blanket of the
venerable Euler-Frobenius-Poincar\'e formal power series approach,
which culminated recently in the developing of symbolic algebra
tools like Maple's {\tt DEtools} and {\tt gfun} and Mathematica's
{\tt HolonomicFunctions}. Such programs easily identify the
general formal series solutions of our special ``Okubo-type"
equations -- see  Section \ref{s:Oku}.

 \vskip.1truecm{\bf Non-dominant/minimal solutions}. There is however one corner not covered well by the above mentioned
 blanket!  Our generating function must be  in the ergodic case the unique (up to a constant) analytic
solution in the unit circle having positive Taylor coefficients around $0$ (as ensured by
ergodic theory), despite the fact that the equation under study has an irregular singularity at
$0$. This phenomenon  seems  less studied,
 and thus, while  computer programs will provide the
subspace of all (typically divergent) power series at $0$, they
will not provide direct information on the unique  generating
function/analytical  solution  with positive coefficients at $0$
that is of interest in probability and combinatorics.

 This difficulty  has already been investigated by \citet{pearce1985holomorphic} in a remarkable uncited paper, where he notes that at an irregular singular point there {might} be two  ``formal power
series solutions that are divergent (have zero radius of
convergence) but that nevertheless may be combined into a
holomorphic solution with a positive radius of convergence in the
vicinity of the singularity. This phenomenon does not appear to
have been noted previously in the literature. The approach of
expressing an arbitrary solution about an irregular singular point
as a sum of solutions of known asymptotic behaviour is not very
convenient for the investigation of regular solutions. In
practice, the determination of regular solutions  centers about
the indicial equation  and around the determination of the
non-dominant/minimal solution for the recurrence satisfied by the
formal series coefficients, for example see \citet[Ch
3]{pearce1985holomorphic}.

\vskip.1truecm{\bf The continued fraction construction of minimal solutions}. For  second order recurrences, the  minimal solution
 may be obtained  by a theorem of Pincherle,
\citet[Thm B.4, pg. 403]{jones1980continued}, or
\citet{gautschi1967computational}, which constructs it by using
continued fractions.

Phung-Duc et al.~\citeyearpar{phung2009M,phung2010state} showed  that  with three or four servers, the recurrence satisfied by the  stationary distribution is of second order, and indeed it is a minimal solution of this  recurrence. Furthermore, they  constructed it  as a linear combination of two  explicit independent solutions of the  recurrence   which diverge (tend to infinity as the number of customers in the orbit tends to infinity), following an approach initiated by \citet{choi1998m}.

Beyond second order recurrences, one needs a generalization of
Pincherle's result to higher order recurrence relations, for
example see \citet{pearce1989extended}, and such  a {\bf matrix
continued fractions generalization} was provided by
\citet{levrie1996matrix}.

    {\bf The matrix continued fractions
  approach to retrial queues} (building  on \citet[Thm. 1]{pearce1985holomorphic},
\citet{pearce1989extended}, and
   \citet{levrie1996matrix}) was   introduced by \citet{hanschke1999matrix}, who emphasized  matrix-algorithmic
  aspects, similar to those of computing
  the $R_n$ matrices of Neuts (the blocks of the RG-factorization in \citet{li2002RG}).

  Recently, \citet{phung2010simple,phung2011matrix} and \citet{baumann2010numerical} have provided impressive applications of the matrix-continued approach to retrial queues and to general level dependent QBD's, respectively. These papers turn  convergence results for matrix continued fractions into   a practical tool (quite related to the $R$ and $G$ matrix algorithmic approach of Neuts).

  %Note also the absence even for $\c=1, 2$ of formulas %for  {the} waiting time and busy period distributions

Many other intriguing facts about retrial queues-- exact results,  asymptotics, approximations and  probabilistic limit theorems -- may be found in \citet{falin1997retrial,artalejo2002standard}, and in the more recent literature.
    Let us end by {listing} a few of these intriguing facts:
 \be \im
 {In the limit to zero retrial rate, \citet{cohen1957basic} discovered} that certain retrial queues  without  a waiting buffer behave like a classic Erlang-loss
system. However, the retrial buffer does not simply disappear, but
gives rise to an increase in the  arrival rate to the system (for
a fine analysis of this phenomenon in the Halfin-Whitt regime, see
\citet{avram2012loss}, \citet{janssen2012asymptotic})).

Note that the Cohen discontinuity from \citet{cohen1957basic} when
the retrial rate converges to $0$ has only been established for
exponential transition times. The matrix approach advocated below
might be useful for investigating this phenomenon under phase-type
transition times.

\im A fundamental result with one server is the {``}stochastic
decomposition" of the generating function as a product, which
contains the generating function of a limit model, for example see
\citet{artalejo1994stochastic,atencia2003queueing,
krishnamoorthy2012queues}).

\im Especially interesting is the  refining of   the existing simple
 approximations like the RTA,
 or constant retrial rate  approximation, and the  Fredericks and Reisner approximation -- see \citet[Ch 7]{wolff1989stochastic},
    and the \citet{grier1997time} approximation.% -- see section \ref{s:MW}.

\im Finally, the role of censoring, which turned out to be crucial
in the analysis of certain retrial models such as
\citet{liu2010analyzing,liutail}, and \citet{kim2012tail},
deserves further investigation.

\im The relation of the matrix continued approach  to the Laguerre series expansions  of \citet{cohen1957basic}
 is  not understood.

  \ee

\vskip.1truecm{\bf Our motivation}: The results mentioned above
have typically been investigated for  ``simple" particular cases
of retrial models. We feel  that this sometimes obscures the
underlying mathematical structure, and therefore decided to start
again at the base:
 Kolmogorov's linear differential
systems for the generating function of   stationary probabilities
for QBD's \eqref{linsys} and retrial queues \eqref{V},
respectively. As a first bonus, we compute the dominant
singularity.

As a second bonus, we hope to  obtain in the future  a full
asymptotic expansion of the stationary probabilities, extending
results  of~\citet{liutail}, and \citet{kim2012tail}. In
principle, this could be achieved by applying   the    vector
Frobenius method, a new
  treatment that has been provided recently by \citet[Theorem 1]{kim2012tail},  but this requires further work, as
  explained below.

\vskip.1truecm
 {\bf Contents}: Retrial queues are a particular
case of  quasi-birth-and-death processes (QBD's),  and some
background material on this family is presented in Section
\ref{s:QBD}.

Our key definition, a general  ``advanced multiserver retrial
model", is presented in Section \ref{s:ret}. The generality
adopted here helps us to organize various open  questions
 (the study of
which we hope to undertake in the future). Section~\ref{s4} deals
with the special Markovian case, a  linear   level-dependent QBD
process, which subsumes several
 interesting queueing features like geometric acceptance, abandon,
 feedback etc.

 Section \ref{s:dif} collects some simple
but fundamental features of our differential system, namely:
 \be
 \im the matrices intervening in the Kolmogorov equations --- see Lemmas~5.1, 5.2,  %\ref{l:d0}, \ref{l:d1}
 \im  a discussion of the corresponding singularities, and
 \im  the stability condition ensuring analicity of the generating function in the unit circle --- see Remark~6.2. % \ref{r:st}.
 \ee

 Section~\ref{s:per} contains results for persistent models.

Section \ref{s:Oku} makes an interesting observation that
persistent retrials with no feedback lead to a recently introduced
type of ``Okubo-type hypergeometric systems" (which always have an
irregular singularity at $0$ when $\c \geq 3$).

Section \ref{s:asy}   considers the  asymptotic behavior of
stationary probabilities (assuming a positive retrial rate). Note
that the cases $s=1,  s=\I$,  {and} $1<s< \I$ correspond
respectively to one variable scalar generating functions
(extensively-studied --- see for example
\citet{odlyzko1995asymptotic,flajolet2009analytic}), two variable
scalar generating functions, essentially an open problem, whose
study has only recently been tackled by
\citet{baryshnikov2004convolutions,pemantle2005twenty,
 pemantle2010analytic,raichev2010asymptotics}, and
\citet{baryshnikov2011asymptotics},
 and to one variable vector generating functions.
 In this latter case,
asymptotic expansions are available in principle by the    vector
Frobenius method
 ---
see \citet{coddington1955theory,wasow2002asymptotic}, and
\citet{barkatou1999algorithm}. First terms of such expansions
have already been obtained  by~
\citet{liu2010analyzing,liutail,kim2012tail}, and we are working
currently on extending these to full asymptotic expansions.

%The two servers case is reviewed in Section %\ref{s:two}.

\section{ {Quasi-birth-and-death} processes: a general framework for bivariate Markovian modelling  \la{s:QBD}}

{\bf  {Introduction:}} Many important stochastic models involve
multidimensional random walks with discrete state space, whose
coordinates split naturally  {into:}
 \be
 \im {\bf one infinite valued coordinate} $N(t) \in {\N =\{0, 1, 2, \ldots \}}$, called \textit{level}, and
 \im  ``the rest of the information" $I(t)$, called \textit{phase  {or background}}, which typically takes a {\bf finite number of possible values}.
 \ee

 Partitioned according to the level, the
infinitesimal generator $Q$ of such a Markov process
$(I(t),N(t))$ is a {\bf block tridiagonal} matrix, called
 \textit{level-dependent} QBD
 generator (LD-QBD):
\begin{equation} \label{QBD}
    Q = \left [ \begin{array}{ccccccc}
\B_0 & \A_0 \\
\C_1 & \B_1 & \A_1 \\
& \C_2 & \B_2 & \A_2 \\
& & \ddots & \ddots & \ddots
\end{array} \right ]
\end{equation}
(recall that $Q$ is a matrix with nonnegative off-diagonal
``rates/weights'', and with row sums equal to $0$). {In the
above and throughout the paper, all unspecified entries of a
matrix are zero.}

{ LD-QBD processes} share the {``skip free"}
structure of {birth-and-death-processes;} however, the
``{\bf{level transition weights}}" $\A_\j,\B_\j,\C_\j$
are now {\bf matrices}, inviting one to enter the noncommutative
world.

{\bf Stationary  {distribution:}} One challenging problem of
great interest
 for {\bf positive recurrent} {LD-QBD} processes is the determination of
the stationary distribution $\vc{\pi} = (\vc{\pi}_0, \vc{\pi}_1,
\vc{\pi}_2, \ldots)$  partitioned by level, where {$\vc{\pi}_{\j}
= (\pi_{0,\j}, \pi_{1,\j}, \pi_{2,\j}, \ldots )$.} The equilibrium
equations
\begin{equation}\label{piq}
\vc{\pi} {Q} = 0
\end{equation}
 in partitioned form yield the second
degree vector recursion:
\begin{equation}\label{genkol}
\begin{cases} \vc{\pi}_{{\j}-1} \A_{{\j}-1} + \vc{\pi}_{\j} \B_{\j} +
\vc{\pi}_{{\j}+1} \C_{{\j}+1} = 0, \qquad {\j} =  1,2, {\ldots} \\
 \vc{\pi}_{0} \B_{0} +
\vc{\pi}_{1} \C_{1} = 0\end{cases}
\end{equation}
% where $\mathbf{\pi}_{-1}$ is a vector of $0$'s.
{Note the absence of an initial value for $\vc{\pi}_{0}$!
(Instead, the sequence is determined by the requirement of
summability  of $\vc{\pi}_{\j}(i)$ in $\j, \forall i.$)}
\medskip

{{\bf LD-QBD  {processes:}}} Level independent QBD processes have
been intensively studied, but the {level-dependent} case is less
understood due to the generality of the model. For some recent
works discussing LD-QBD processes, see
\citet{bright1995calculating,ramaswami1996matrix,
li2004two,li2006BMAP,li2010constructive}, and
\citet{avram2011symbolic}.

Despite recent theoretic and algorithmic progress in the study of
formal series solutions, there is an intriguing scarcity of
analytic works in the particular ``holonomic" case of linear or
polynomial dependence on the level, although such recurrences and
the corresponding generating functions have been central in
classic works of Euler, Gauss, Riemann, Liouville, Klein, Fuchs,
Frobenius, Pincherle, Poincar\'e, Lie, Hilbert, Birkhoff, etc (see
\citet{gray2008linear} for a delightful account), and despite
recent theoretic and algorithmic progress in the study of formal
series solutions. Narrowing the gap between this venerable
mathematics and our intriguing applications has in fact been our
motivation for this work.

We end by noting that  very efficient numerical approaches for computing stationary probabilities  for
 {level-dependent} QBD's have been provided recently by
\citet{baumann2010numerical,phung2010simple}, and
\citet{phung2011matrix}.

\section{The {GI/G/\c/K queue} with  retrial, Bernoulli acceptance,  abandon
 and feedback
 \la{s:ret}}

The inclusion of this general model in a  paper which deals
with a modified $M/M/\c/K$  Markovian model requires some
justification.
We are motivated by the fact, as recognized in Hanschke (1999), that when the primary area has no free servers, first time arrivals face one of three choices:
a)
wait in  the primary area (in the priority queue), b) abandon, and
c) join the orbit for later retrials. The same three choices may
be made after each service {by the customer who has just
completed his service (since feedback is allowed in our model),
and  retrial customers  also
face the first two choices.} However, several recent fundamental
papers still assume that one of these three options (which changes
from paper to paper) may not happen.

One such example, and for us an especially important one, is the generalized Pollaczek-Khinchine formula of Atencia and Moreno~\citeyearpar{atencia2003queueing}. We feel that their generating function formula \citet[Thm 3.1]{atencia2003queueing}
marks
 an important moment in the evolution of retrial queues, since a) it is general enough to include many particular cases studied previously, and b)  it allows for general service times,  so that it may be applied directly to empirical data.

 Unfortunately, some of the three natural retrial options are still assumed not to happen in  \citeyearpar{atencia2003queueing};
this motivated us to formulate
the   general  three-choice model below, so that  we may state an important open problem:
\beQ Generalize the Atencia-Moreno formula for the one server, and eventually for the multi-server case of the model below. \eeQ

\begin{Def} %\label{o-r}
 The  {\bf overflow GI/G/\c/K  queue with retrial rate $\n$, orbit abandon probability $p$, Bernoulli acceptance probability $p_a$ and
feedback to orbit probability $\th$} consists of two facilities: (1) a
``primary service area'' with $\c$ service places (or servers) and $K$  extra waiting spaces (referred to as a priority queue), the only cases studied here being $K=0$ and $K=\infty$.   ``Primary"
customers arrive to the primary service area according to a renewal process with rate $\l $ (the
case of interest for us being Poisson arrivals); and  (2) an additional
infinite ``overflow buffer" (referred to as an orbit).
\end{Def}

 The main {``}activities" of interest in the system are:
 \begin{description}
   \item[(1)] {\bf Arrivals} to the primary  area  with rate $\l $.
   If at least one server is
 free  upon arrival, {the customer is} accepted  with
probability $ p_a$ for first time customers,  sent to orbit with
 probability $\T  p_a,$ or forced to leave the system for some obscure reason with
 probability $\bar p_a=1-p_a-\T  p_a$.
 However,
if  all servers are busy  {upon arrival, but $0<K<\I$ and there is still space in the priority queue, then these
probabilities become $ \b_0 ,\T \b_0, \bar \b_0 =1 -\b_0 -\T \b_0
$ (acceptance means in this case that customers will wait in the
priority queue of the primary area). To simplify notation, we will assimilate this case with  the first, i.e. assume that $p_a=\b_0, \T p_a= \T \b_0$, \ldots (since  we consider  only  $K=0$ below), keeping in mind that in certain applications a distinction might need to be made.} Finally, if the servers and the priority queue are blocked,  these
probabilities become $ \a_0=0 ,\T \a_0, \bar \a_0 =1 -\T \a_0
$. We   denote the rates of the
three possibilities, with available servers and after blocking  {of the
primary area}, respectively, by
\[
\begin{cases}\l_q= \l p_a&\text{join the primary area}\\ \l_a= \l
\bar  p_a&\text{leave the system} \\ \l_o= \l \T  p_a &\text{join
the orbit}\end{cases}
\q \begin{cases}\l_{qb}= 0\\ \l_{ab}=\l \bar  \a_0 \\
\l_{ob}= \l \T \a_0 \end{cases}
\]

\end{description}

\begin{remark} Our model is related to a model proposed by
\citet{hanschke1999matrix}, essentially by adding feedback, and to the model of \citet{atencia2003queueing}.
To clarify this, we include here first Hanschke's, and then Atencia and Moreno's notations for the parameters
(they correspond to  certain probabilistic interpretations, which do not play a role in the mathematical analysis).

\small{
\[
\hspace*{-5mm}
\begin{cases}\l_q= \l p \; ({\mathcal H})=\l \;
({\mathcal AT})&\text{join the primary area}\\ \l_a= \l (1-
p)(1-r_1)\;
({\mathcal H})=0 \; ({\mathcal AT})&\text{leave the system} \\
\l_o= \l (1-  p)r_1 \; ({\mathcal H}) =0 \; ({\mathcal
AT})&\text{join the orbit}\end{cases} \q \begin{cases}\l_{qb}= 0
\; ({\mathcal H})=\l q \; ({\mathcal AT})\\ \l_{ab}=\l (1-q_1) \;
({\mathcal H})=0 \; ({\mathcal AT})\\ \l_{ob}= \l q_1 \;
({\mathcal H})=\l p \; ({\mathcal AT})\end{cases}
\]
}

\noindent Clearly, adopting a unified notation would be beneficial.
 \end{remark}

\begin{description}
\item[(2)] {\bf Leaving the orbit for retrial or abandon}. \ The
orbit may decrease by one at an ``affine" rate $\n(i)= \n^{(0)}
\1_{i
>0}+ i \n$. The first term, which may be interpreted as the effect
of service by one  secondary server, will be assumed to be $0$
throughout
 the paper. The second term is due to each
customer  reevaluating his position of  remaining in the orbit  at rate $\n$.  {At these evaluation times, he is informed first on the state of the primary area (free service, free waiting space, blocked primary area). Then, he abandons the system with  probability $ \bar \a$ or $ \bar p, $  depending on whether   the primary area is blocked or not. In the second case the remaining proportion of $  p=1-\bar p$ will go to the primary area. In the first, the remaining  proportion of $  \T \a=1-\bar \a$ will remain in the orbit. All customers in the orbit repeatedly reevaluate their position  until a free place is secured in the primary area (either an idle server or in
the priority queue), referred to as a successful retrial, or until abandon. }  We will denote the respective rates of going to the primary area and abandoning
before and after blocking by
\[
    \begin{cases}\n_q= \n p &\text{join the primary area},\\ \n_a= \n \bar p &\text{leave the system}\end{cases} \begin{cases} \n_{qb} = 0 \\
\n_{ab}=\n \bar \alpha \end{cases}
\]

 {When $\bar \alpha=\bar p=0,$ the retrials will be called persistent}.

In \citet{hanschke1999matrix}, the respective rates are denoted by
\[
    \begin{cases}\n_q= \n p&\text{join the primary area}\\ \n_a=
\n (1-  p)(1-r_2)&\text{leave the system} \end{cases} \q
\begin{cases}\n_{qb}= 0\\ \n_{ab}=\n (1-q_2) \end{cases}
\]
\end{description}

\begin{remark}
 Note that % {upon a successful retrial,}
 we have not included a third
possibility of {rejoining the orbit}, which could be appropriate in a discrete time model but not
 in a continuous time model, since staying put does not produce a
transition rate in this case. For arrivals, it is also possible to similarly omit
the abandoning arrivals in the specification of the model. This
amounts mathematically to assuming that $\bar p_a=1$, and letting
$\l$ denote the rate of nonrejected arrivals, as   we will do from
now on.

 \end{remark}

\begin{remark}
The parameter $\n$ may be interpreted as the total
activity rate per individual in the orbit (retrial + abandon due
to impatience), and  {$p$ and} $\bar{p}$ as the respective
probabilities of these activities; normally, one would expect
$\bar{\a} \geq \bar{p}$  {i.e.}  a bigger probability of
leaving due to congestion.

Note also that {from the mathematical modelling point of
view,} impatience is identical to having an additional infinite server
dispatching the customers.

\end{remark}

\begin{description}
\item[(3)] {\bf  Service and feedback}. For each of the accepted
first time or  retrial customers, the service time  follows a
general distribution with affine service rate (=inverse of
expected service time) of rate $\mu(i) =\mu^{(0)} \1_{i>0} + i
\mu$. We will  consider mostly the case when $\mu(i)  $ is either
constant or linear,  with exponential service times.
  All
service times are independent, and are also independent of the arrival,
retrial and routing  processes.

{\bf Feedback: leaving the system/joining the orbit after
service}.  After service in the primary  area,   {the
customer may leave} the system forever with probability  $\bar
\theta$, request an  additional service in the primary  area
 with probability $ \theta$,     or
 join  the orbit with probability
$\T \theta=  1 -\theta -\bar
\theta$ (the ``no feedback" factor $\bar \theta$ is $1$ in the standard model).
We will  denote the rates of  these three possibilities by
$$\begin{cases} \m_q= \m
  \theta\\\m_a= \m \bar  \theta\\\m_o= \m  \T \theta
\end{cases}$$
  Here we have followed the notation of \citet{atencia2003queueing}, while adding the third case with parameter $\T \theta \geq 0$ (note that since feedback is not affected by possible blocking of the primary
area,  we only need one set of rates in this case).
%and assume for simplicity below that $\T  \theta=0$.

\end{description}

%\end{Def}

We included the above general model because we hope that it will encourage the emergence of generalizations of the Atencia-Moreno formula \citeyearpar[Thm 3.1]{atencia2003queueing}. However, as is common in the multiserver literature, we will now go on to consider the M/M/s/K retrial model.

%From now on we will consider,  like  most of the multiserver
%literature,   the M/M/s/K retrial model  with exponential
%arrivals, services and retrials,  with $ K=0$ or $K=\I$, and with
%$\m(i)$ and $\n(i)$
%either  linear or constant.
%We included however the general model because we hope that this would
%encourage the emergence of generalizations of the Atencia-Moreno formula
%\citeyearpar[Thm 3.1]{atencia2003queueing}. %, at least under the assumption of phase-type %transition times.

\section{The Markovian model with exponential transition times}
\la{s4}

We now consider the case of exponential arrivals, services and retrials.

\begin{Def}
An {\bf affine death QBD process} is a QBD process with affine
  dependence on the level:
 \beq \la{linQBD}
 && \A_{\j} =\A, \q \C_\j= \j \C + \C^{(0)} 1_{\{\j >0\}} , \\&& \no \B_\j=\B  - \T \A  -\j \T \C-\T \C^{(0)} 1_{\{\j
>0\}}
 \eeq
 where $B$ is a conservative generator, and $\T \A$, $\T \C^{(0)}$ and $\T \C$ are diagonal matrices containing on their
diagonals  the sum of the rows  of $\A$, $\C^{(0)}$ and $\C$,
respectively.
\end{Def}

The two terms in $\C_\j$ correspond to the cases of
``{\bf independent retrials}" and ``{\bf one retrial dispatcher}"
in
the %secondary service area/
orbit, respectively, and we here are mainly interested in the
linear retrial case with {$\C_\j= \j \C$ and $\C^{(0)}=0$,} and in
the
 constant  retrial case $\C=0$, when  we will denote  $\C^{(0)}$  simply
by $\C$.

Let  {$(N(t),I(t))$} denote respectively the numbers of
customers in the  {orbit and in the service area} at time $t$
  for a
Markovian retrial system. Consider first the case {with $s$
servers and no extra waiting spaces ($K=0$)}, assuming for now
$\C^{(0)}=0$. Then, the retrial process {$(N(t),I(t)) $} with
 {respective geometric losses $\bar{p}_0 = 1-p \in [0,1]$ and
$\bar{\alpha}_0 =1- \a_0 \in [0, 1]$, orbit abandon probabilities
$\bar{p}=1-p \in [0, 1]$ and  $\bar{\alpha}=1-\alpha \in [0, 1]$,}
 and  feedback $\th \in [0, 1]$
 is a linear death {QBD process with departure matrix}
 $\C_\j = \j \C$,
 where:
{\[ {\C} = \left[\begin{array}{ccccc}
%\begin{equati\ddotson*}
%\begin{pmatrix}
   \n   \bar p    &   \n   p &       \\
    &      \n   \bar p     &    \n   p &   \\
    &    &   \ddots    &   \ddots  &    \\
    &   &     &   \n   \bar p   &     \n   p  \\
    &   &     &   &   \n \bar \a \\
\end{array}\right]
\] }
 and with  {arrival and ``loop"} (transitions  without level
change) matrices:
 \beq \la{DefA} \A= \left[\begin{array}{cccccc}
   \l  \T  p_a    &    \\
   \m   \Tth    &       \l  \T  p_a   &     \\
    &      2  \m   \Tth   &   \l  \T  p_a   &   \\
&  &   \ddots  &   \ddots    &  \\
& & &  (s-1) \m \Tth & \l  \T  p_a  &  \\
     &     &   &  &    \c \m   \Tth     &     {\l \T \alpha_0} \\
\end{array}\right]
 \eeq
 \beq
 {\B} = \left[\begin{array}{ccccccc}
%\begin{equati\ddotson*}
%\begin{pmatrix}
-\l  p_a &  \l  p_a  &  \\
\m  \bar  \theta  & -(\l  p_a + \m  \bar  \theta)  &  \l  p_a &     \\
  & 2 \m  \bar  \theta  & -(\l  p_a + 2\m  \bar  \theta) &  \l  p_a &    \\
  &  &   \ddots  & \ddots &  \ddots   &  \\
  & & &    (s-1) \m  \bar  \theta & -(\l  p_a + (s-1)\m  \bar  \theta) & \l  p_a\\
  & & & &   \c \m  \bar  \theta  & - \c \m  \bar  \theta \\
\end{array}\right]\eeq
where $p+  \bar{p}=1$  {and $\bar p_a+p_a=1$.}

%\begin{remark} {This remark can be deleted.} Note that $B$ is %a
%conservative semigroup generator matrix.\end{remark}

\begin{example} In the  retrials model   with $p_a=1$ (all accepted when
 {a server} is available),
{%\large
the arrivals $\A_\j=\A$ and ``loops" are:
\[
\A= \left[\begin{array}{cccccc}
  0     &   \\
   \m \Tth    &       0   &     \\
&   2  \m   \Tth   &   0     &   \\
& & \ddots  &   \ddots    &  \\
& & & (s-1) \m \Tth & 0 & \\
& & & &  \c \m   \Tth     &     \l_{ob} \\
\end{array}\right]
\]
\[
 {\B} = \left[\begin{array}{ccccccc}
%\begin{equati\ddotson*}
%\begin{pmatrix}
-\l  & \l &       \\
\m  \bar  \theta  & -(\l + \m  \bar  \theta)  & \l &     \\
  & 2 \m  \bar  \theta  & -(\l   + 2\m  \bar  \theta) & \l  &  \\
  & &  \ddots & \ddots & \ddots &  \\
  & & & (s-1) \m \theta  & -(\l   + (s-1)\m  \bar  \theta) & \l \\
        &&& &   \c \m  \bar   \theta  & - \c \m  \bar  \theta \\
\end{array}\right]
\]
}

\end{example}

\begin{remark} Our model
includes three parameters which are less studied: \be
 \im
 the  {``}abandon retrials if blocked" parameter $\bar \a$, % which allows to include the case of
%"geometric orbits",

\im  the  {``}impatience abandon retrials" parameter $\bar p$,
which allows us to also include the easier case of Markov modulated
M/M/$\I$ queues, characterized by $\bar p=1$, which is considered
for example in \citet{economou2005generalized},
  and
\im the feedback parameter $\th$.
\ee\end{remark}

\section{The  differential system for the generating functions %case $K=\c$
\la{s:dif}}

Let $\pi_{i,{\j}}$ denote the stationary probabilities of having
$i$ customers in the primary area and $\j$ customers in {the}
orbit, which satisfy the recursion \eqref{genkol}. A classic
approach for tackling  this is via
 the {\bf generating functions}
 \begin{equation}\label{deffoncgen}
{ p_i(z)=\sum_{{\j}=0}^\infty \pi_{i,{\j}} z^{\j}, {i=0, 1, 2
\ldots}  \; \Longleftrightarrow \; }
\mathbf{p}(z)=\sum_{{\j}=0}^\infty \vc{\pi}_{{\j}} z^{\j},
\end{equation}
where $\mathbf{p}(z) =  {(p_0(z), p_{1}(z), \ldots)}$.

 \begin{lemma}
 \label{l:d0}  a) For the affine death QBD process in \eqref{linQBD}  {(with finite $s$ and finite $K$)}, the recursion
  \beq\label{linrec}
\begin{cases} \vc{\pi}_{{\j}-1} \A + \vc{\pi}_{\j} (\B - \T A -{\j} \T C
-\C^{(0)}) +
 \vc{\pi}_{{\j}+1} (({{\j}+1})\C+ \C^{(0)}) = 0, \; {\j} =  1,2, \ldots\\
 \vc{\pi}_{0} (\B - \T A) +
\vc{\pi}_{1} (\C + \C^{(0)})= 0\end{cases}
\eeq
yields  a  linear differential system
\begin{eqnarray}&&
\la{linsys}
   {\mathbf{p'}(z) V(z) = \mathbf{p}(z) U(z) +  \vc{\pi}_{0} (\T \C^{(0)}-  z^{-1} \C^{(0)})},
\end{eqnarray}
where $V(z)$  {and} $U(z)$ are square matrices
 of order $\c +1+ K$, given by:
\beq \la{VU} &&  {U(z)= B + z \A - \T \A + z^{-1} \C^{(0)}  - \T
\C^{(0)}, \q V(z) = z \T \C -  \C} \eeq

b) In the particular case of the    model of Definition~1 % \ref{o-r}
with exponential transition times, $ \C^{(0)}=\T
\C^{(0)}=0$, $K=0,$ and $1 \leq \c < \I,$ (\ref{linsys})  holds
with $V(z)$ and $U(z)$  given by:
 \be \im For $\c =1,$ \beq &&V(z)=\n
\left[\begin{array}{cc}
%\begin{equati\ddotson*}
%\begin{pmatrix}
   z - \bar p  &  -p    \\
   0     &     {\bar{\a}  (z-1)} \\
\end{array}\right], \q {U(z)} =
\n \left[\begin{array}{cc}
%\begin{equati\ddotson*}
%\begin{pmatrix}
\Tl( z \T  p_0-1)&     \Tl  p_0             \\
       \Tmu (\bar  \th +\Tth  z)
       & \Tl_{ob} (z-1) - \Tmu(\bar  \th +\Tth)
\end{array}\right]. \no \eeq
\im {For $2 \leq \c < \infty$:
\beq &&{V(z)=\n^{-1}(z \T \C -  \C)}=
\left[\begin{array}{ccccc}
%\begin{equati\ddotson*}
%\begin{pmatrix}
   z - \bar p  &  -p  &     \\
       &      z - \bar p   &   -p   &   \\
&  &   \ddots  &   \ddots    &   \\
& &  &   z - \bar p&   -p  \\
   &&&  &    {\bar{\a}  (z-1)} \\
\end{array}\right] \la{V}
\eeq } \ee

 {\small
\begin{align*}
& {U(z)} %- z^{-1} \C^{(0)}  +  \T \C^{(0)}
= \n^{-1}(\B- \T \A + z \A) = \\
& \left[\begin{array}{cccccc}
%\begin{equati\ddotson*}
%\begin{pmatrix}
\Theta_0 & \Tl  p_0  & \\
\Tmu (\bar  \th +\Tth  z)  & \Theta_1 & \Tl  p_0  &  \\
 & 2 \Tmu   (\bar  \th +\Tth  z)  & \Theta_2 &  \Tl p_0 &  \\
& & \ddots & \ddots & \ddots & \\
&&&  (\c-1) \Tmu (\bar  \th +\Tth  z) &  \Theta_{s-1}  &  \Tl p_0  \\
&&&&  \c \Tmu (\bar  \th +\Tth  z) & \Tl_{ob} (z-1) -\c \Tmu (\bar  \th +\Tth)
\end{array}\right] \no
\end{align*} }
 where
 $\Tl=\fr{\l}{\n  }$, $\Tmu=\fr{\m }{\n }$ and \beq \la{eq:Th} \Theta_k=\Tl( z \T  p_a-1) - k\Tmu (\bar  \th +\Tth), \; k=0, 1, \ldots, (s-1).\eeq

\medskip

The bivariate {\bf generating function}
 \begin{eqnarray*}
\vp(y ,z):=  {\sum_{i=0}^{\c} \sum_{\j=0}^{\infty}}
\pi_{i,{\j}} y ^i z^{\j} \end{eqnarray*} satisfies the equation: \beq
\label{e:bigen}
&& \n(z-\bar p - p y) \vp'_z(y ,z)+\n y^\c (\bar p- \bar \a+ z(\bar \a-1)-p y ) p_\c'(z) = \\
&&\hspace*{1cm} \l \vp(y ,z)(\T p_0 z  -1+ p_0 y)+\mu
 \vp'_y(y ,z)(\bar \th + \Tth z -y(\bar  \th +\Tth)) \no \\
&&\hspace*{1.5cm}  + y^\c(z(\l_{ob}-\l \T p_0)-(\l_{ob}-\l)-\l p_0 y) p_\c(z) \no\eeq

c) In the infinite buffer case $K=\I$, the differential system $
  \mathbf{p'}(z) V(z) = \mathbf{p}(z) U(z)$ formed from {\eqref{e:int}, \eqref{e:ext} and \eqref{e:bd}} involves now infinite square matrices.  For $\c=1$ for example,
the  $V$ and $U$ matrices
 are  {given respectively by}:
\begin{eqnarray*} &&V(z)=
\left[\begin{array}{ccccc}
%\begin{equati\ddotson*}
%\begin{pmatrix}
   z - \bar p  &  -p   & \\
   &     ({\bar{\a} +\a) z-\bar{\a}}&-\a&  \\
   & &     ({\bar{\a} +\a) z-\bar{\a}}&-\a& \\
   & & & \ddots&\ddots \\
   &&&& \ddots
\end{array}\right]
\end{eqnarray*}
  {and}
 { \begin{eqnarray*} {U(z)} = \left[\begin{array}{cccccc}
%\begin{equati\ddotson*}
%\begin{pmatrix}
   \Tl( z \T  p_0-1)&     \Tl  p_0  & \\
   \Tmu (\bar  \th +\Tth  z) & \Tl (\T  \a_0 z-\T  \a_0-\a_0 ) - \Tmu(\bar  \th +\Tth)&\Tl \a_0&  \\
   &\Tmu (\bar  \th +\Tth  z)    &  \Tl (\T  \a_0 z-\T  \a_0-\a_0 ) - 2\Tmu(\bar  \th +\Tth) &\Tl \a_0& \\
   &   & \ddots&\ddots&\ddots
\end{array}\right] \no \la{V1}\end{eqnarray*}
}

\begin{comment}
When $K=\I,$ the bivariate {\bf generating function}
 \begin{equation}\label{e:bigen}
\vp(y ,z):=\sum_{i,{\j}=0}^\infty \pi_{i,{\j}} y ^i z^{\j}
\end{equation}
 satisfies the equation: ...
 \end{comment}
\end{lemma}

{\bf Proof:} a) {Multiplying
 the equilibrium equations in  \eqref{linrec} by
$z^\j$ and summing up gives rise to a linear first order
differential system:
 \begin{eqnarray*}&&
\sum_{\j=0}^\I z^\j \vc{\pi}_{{\j}-1} \A + z^\j \vc{\pi}_{\j} (\B-\T
\A) -  \vc{\pi}_{\j} z^\j  (\j \T \C+   \T \C^{(0)} 1_{\{\j >0\}}) +
 z^\j \vc{\pi}_{{\j}+1} (({\j}+1) \C +  \C^{(0)}) = \no \\
 && z
  \vc{p}(z) \A + \vc{p}(z) (B - \T \A-  \T \C^{(0)}) +
  \vc{\pi}_{0} \T \C^{(0)} - z  \vc{p}'(z) \T \C+ \vc{p}'(z)
  \C + z^{-1} (\vc{p}(z) -\vc{\pi}_{0})\C^{(0)} = 0
\end{eqnarray*}
 {which yields \eqref{linsys}.} }

b,c) It is convenient to associate a number with each of the
activities of the system. Let us denote by 1, 2, \ldots, 6,  the
inflow activities associated  {with} the rates $\l_w,
\l_{ob},\mu_o,\mu_a,\n_a,\n_w$, respectively, and by {$\langle
12 \rangle$, $\langle 34 \rangle$ and $\langle 56 \rangle$} the
corresponding outflow activities (the distinction between the
blocked and non-blocked cases is not made here).

 It is easy to check that the {\bf generating function terms associated  {with} respective activities} are
$$\begin{array}{|c|c|c|ccc|} \hline 1&\l p_0 p_{i-1}(z)&\l \a_0 p_{i-1}(z)&i \geq &\c +&1
\\ 2 &\l \T  p_0 z p_{i}(z)&\l \T  \a_0 z p_{i}(z)&i \geq &\c&\\  {\langle 12 \rangle} &\l (p_0+\T  p_0 )  p_{i}(z)&\l (\a_0+ \T  \a_0)  p_{i}(z)&i \geq &\c&\\\hline 3 &(i+1) \mu  \Tth z p_{i+1}(z)& \dots&i \geq& 0&
\\4&(i+1) \mu \bar \th p_{i+1}(z)& \dots&i \geq & 0&\\
 {\langle 34 \rangle} &\mu i(\bar  \th +\Tth )  p_{i}(z)&\mu i
(\bar \th +\Tth ) p_{i}(z)&i \geq &\c&\\\hline 5&\n \bar p
p_{i}'(z)&\n \bar \a p_{i}'(z)& i \geq &\c&\\6&\n  \a  p_{i-1}'(z)
&\n  \a  p_{i-1}'(z)&i \geq &\c +&1  \\
 {\langle 56 \rangle} &\n (p+\bar  p) z p_{i}'(z)&\n (\a+\bar  \a) z p_{i}'(z)&i \geq &\c& \\
\hline
\end{array}$$

The equilibrium equations may be written  symbolically as \beq
\la{e:sy} -6 +  {\langle 56 \rangle -5=1+2 - \langle 12 \rangle
- \langle 34 \rangle} + 4 +3 \eeq

 In the interior,  for $i=1, \ldots, \c-1$, using the table this results in
\beq \la{e:int} &&
\n\le(- p  p_{i-1}'(z)+  ((p+\bar  p) z - \bar p ) p_{i}'(z)\ri)= \\
&& \l p_0 p_{i-1}(z)+  \l(\T  p_0 z - p_0-\T  p_0 ) -\mu i \le(\bar  \th +\Tth  \ri)p_{i}(z)+(i+1)
(\bar  \th +\Tth  z) p_{i+1}(z) \no\eeq
as  stated.

The exterior equilibrium equation for $i=\c+1, \c+2, \ldots$ (also
represented symbolically  {in} \eqref{e:sy}) is, due to the
changes in rates:
 \beq \la{e:ext}
 &&\n\le(- \a  p_{i-1}'(z)+  ((\a+\bar \a) z - \bar \a ) p_{i}'(z)\ri)= \\
 && \l \a_0 p_{i-1}(z)+  \le(\l(\T  \a_0 z - \a_0-\T  \a_0) -\mu i (\bar  \th +\Tth ) \ri)p_{i}(z)+(i+1)  (\bar  \th +\Tth  z) p_{i+1}(z)
 \no\eeq

Finally, the boundary equilibrium
equation for $i=\c$ is
\beq \la{e:bd} &&\n\le(- \a  p_{i-1}'(z)+  ((\a+\bar \a) z - \bar \a ) p_{i}'(z)\ri)= \\
&&\l p_0 p_{i-1}(z)+  \le(\l(\T  \a_0 z - \a_0-
\T  \a_0) -\mu i (\bar  \th +\Tth ) \ri)p_{i}(z)+(i+1)  (\bar  \th +\Tth  z) p_{i+1}(z)\no\eeq

\begin{comment}
\begin{remark} \la{Edefl} The QBD defining blocks may be written as:
\begin{eqnarray*} &&\A=\l_{ob} M +\l \bar{p}  (I-M)+ \m  \theta E_-, \q  \B = \m
\theta (E_- -   E_- T_+) + \l p (T_+ - (I-M)), \\&&\C=  \n_r   T_+ +
 \n_a   (I-M) + \n_{ba} M, \q \C^{(0)}= \s p  T_+ + \s \bar p (I-M) +
\s   \bar{\a} M \end{eqnarray*} in terms of the "basic" matrices defined in
Remark \ref{Edef}.
\end{remark}

Note that  the phase generating matrix $ T_{\j}:=\A_{\j}+ \B_{\j} +
\C_{\j}=(\j \n   + \l p)(T_+ - (I-M))   + \m  (E_- - E_- T_+)$ has
indeed sum of rows $0,$ as it should.

{\bf For the perturbation around $p=0,$ the only difference is that
the orbit is now a Markov modulated M/M/$\I$ queue (with two phases,
corresponding to a full and non-full primary area).

Another possibility is to perturb around a constant retrial queue
with parameters provided by Cohen's formula.}

\end{comment}

\begin{remark} The bivariate generating function equation \eqref{e:bigen} generalizes
the equation \citet[(2.21)]{falin1997retrial}:
 \begin{eqnarray*} && \n(z- y) \vp'_z(y ,z)-\n y^\c (z -y) p_\c'(z)=\l \vp(y ,z)(  -1+  y)+\mu
 \vp'_y(y ,z)(1 -y) +\l  y^\c (z-y) \vp_\c(z).\end{eqnarray*}
 \end{remark}

\begin{remark} {\bf The Markov modulated M/M/$\I$ case.} In the particular case of
a Markov modulated M/M/$\I$ queue, with environment transitions
specified by $\B$,   the arrival and departure matrices $\A$ and $\C$  are diagonal,
$\C^{(0)}=0$, and the differential system \eqref{linsys} for the
generating function
 becomes
\beq \la{eq:mod} (1-z) \mathbf{p'}(z)  C =
 \mathbf{p}(z) ( (1-z) A-  B), \eeq
 or
 \beq \la{gfunMM} \mathbf{p'}(z)   = \mathbf{p}(z) ( \hat {A}-
(1-z)^{-1} \hat {B}),  \q  \hat{A} = A C^{-1}, \hat{B} = B C^{-1},
\eeq
when  $C$  is invertible.

In this case, \citet[Thm 3.1]{o1986m}  expand \eqref{eq:mod}
in power series around $z=1$,  thus getting a recurrence for the conditional factorial moments:
\beq \la{meansMM}
&&\mathbf{p}^{(k)}(1) (k \C -\B)= k \mathbf{p}^{(k-1)}(1) \A, \q
k=1,2, \ldots \\&&\mathbf{p}(1) \B=0\eeq
By the last equation, the phase stationary probabilities $\mathbf{p}(1)$ are the stationary vector of $B$.

 Note that the factorial moments  yield efficient approximations for the stationary probabilities, simply by  {``}shifting the expansion"  {to around} $z=0$.

This resolved an open problem signalled by  Neuts. Citing \citet{neuts1995matrix},  {page 274:}  ``We note that the
infinite-server M/M/$\I$ queue  in a Markovian environment is
surprisingly resistent to analytic solution \ldots Brute force numerical
solution \ldots seems necessary because of the lack of a mathematically
elegant solution.''
It is intriguing to investigate whether a similar approach based on a series expansion around the regular point $1$ could produce the factorial moments for retrial queues as well.%(or $k \mathbf{m}_k \C \; F =
%\mathbf{m}_{k-1} \A \; F,
% k \geq 1$, where $F$ is a matrix with a first column of $1$'s and $0$
%else).
\end{remark}

\begin{remark}
When $K=0,$ the singularities of the system  \eqref{linsys} are the roots of
$Det(V(z))=\bar \a (z-\bar{p})^{\c} (z-1)$. However, it turns out
 that the singularity at $1$ may be
simplified by replacing the last equation by the sum of the
equations, and  dividing by $z-1$.
\end{remark}

\begin{lemma}[Pearce's Lemma]  %\label{l:d1}
{With $K=0$ and $p +\bar p=1$,
 the system}  \eqref{linsys}
is equivalent to the {\bf simplified system} $ \mathbf{p'}(z) V(z) =
\mathbf{p}(z) U(z),$   with {$V(z)$ and $U(z)$} for $\c \geq 2$ given by
%\large
  \beq \la{sysred} &&V(z)=\left[\begin{array}{cccccc}
   z  -\bar p &  -p  &           &    &\cdots &  1   \\
       &      z  -\bar p  &   -p   &  &    &     1   \\
& &   \ddots  &   \ddots    &     & \vdots \\
   && &    z -\bar p& -p &   1  \\
   &&&&  z -\bar p&   1  \\
    &  &&& &     {\bar{\a}  } \\
\end{array}\right]
\eeq
{\small {
   \beq  U(z)=
\left[\begin{array}{cccccc}
\Theta_0 & \Tl  p_0 &   &   &  & \Tl \T  p_0   \\
\Tmu (\bar  \th +\Tth  z)   & \Theta_1 &  \Tl  p_0  &  &  & \Tl \T  p_0 + \Tth \T \mu\\
  &  2 \Tmu  (\bar  \th +\Tth  z) & \Theta_2 & \Tl  p_0  &  & \Tl \T  p_0 + 2 \Tth \T \mu \\
  &   &  \ddots  & \ddots & \ddots &  \vdots \\
  & & & (\c-1) \Tmu (\bar  \th +\Tth  z)  & \Theta_{\c-1} & \Tl \T  p_0 +(\c-1) \Tth \T \mu \\
  &   & &  & \c \Tmu (\bar  \th +\Tth  z) &
        \Tl_{ob} + \c \Tth \T \mu
\end{array}\right ]
\eeq} }

   %(see Remark \ref{simpl}).

   For $\c =1,$
the  {matrices $V(z)$ and $U(z)$
 are given by:}
\beq &&V(z)=\n
\left[\begin{array}{cc}
%\begin{equati\ddotson*}
%\begin{pmatrix}
   z - \bar p  &  1    \\
   0     &     {\bar{\a}  } \\
\end{array}\right], \q {U(z)} =
\n \left[\begin{array}{cc}
%\begin{equati\ddotson*}
%\begin{pmatrix}
\Tl( z \T  p_0-1)&     \Tl  p_0             \\
       \Tmu (\bar  \th +\Tth  z)
       & \Tl_{ob} + \Tth \Tmu
\end{array}\right] \no \la{V1r}\eeq

\end{lemma}

{\bf Proof:}    Let us  replace the last equation  by the sum of the
equations, i.e.
  replace the last columns of {$V$ and $U$} by the sum of the columns.
 Since $\B$ is a conservative generator, this sum equation is
 $\bff p'(z) (z-1) \T \C \bff 1= \bff p(z) (z-1) \T \A \bff 1,$
 and $z-1$ may be simplified,  yielding:
\beq \la{csym} \sum_{i=0}^{\c-1} p_{i}'(z)+ {\bar \a} p_{\c}'(z)=\Tl
 {\Big(} \T  p_0 \sum_{i=0}^{\c-1} p_i (z) +\a_0 p_{\c}(z)  {\Big)}+ \Tth \T \mu
\sum_{i=1}^{\c} i p_i (z)\eeq

 \begin{remark}
The  ``simplification" by $z-1$ of the original system seems to have been first
  noticed in \citet[(4.6)]{pearce1989extended}, who used the first $\c$ equations of \eqref{sysred}
 to obtain
 formulas expressing $p_i(z), i = 1, \ldots,\c$ in terms of $p_0(z),$
 and then used  \eqref{csym} to get  a scalar ODE for $p_0(z)$
  (and  recursion
 for the respective probabilities).
\end{remark}

\begin{remark}
 The system
\eqref{sysred} has been solved in particular cases when its degree is at most two in
terms of classic hypergeometric functions, basically by
``table look up" of the solution
  --- see for example \citet{hanschke1987explicit,falin1997retrial}, and especially  \citet{choi1998m}, who assumes $\bar \a >0$ and includes feedback.

\end{remark}

One must further distinguish the two cases of: \be

\im  ``Non-persistent retrials" with $\bar \a >0,$
{for which} the system is always ergodic. % and  the asymptotic
%behavior is different.

 \im ``Persistent retrials" with
$ \bar \a=0,$ when the determinant is $0$ identically, and a reduction of the dimension by $1$ is possible.

  \ee

  Below, we consider only the second case.
\section{Persistent retrials ($    \bar{\a} =0$) and dimension reduction \la{s:per}}

A special treatment is necessary in the case of ``persistent retrials"  with
$\bar{\a} =0,$ when our system becomes $
\mathbf{p'}(z) V(z) = \mathbf{p}(z) U(z)$ with $V(z)$ given by
  \begin{eqnarray*} &&\left[\begin{array}{cccccc}
   z  -\bar p &  -p  &     &       &     &  1   \\
    &      z  -\bar p  &   -p   &  &   &     1   \\
&  &   \ddots  &   \ddots    &    &    \vdots \\
 & & &   z -\bar p& -p &   1  \\
 & & & &  z -\bar p&  1  \\
 & & & & & 0 \\
\end{array}\right]
   \end{eqnarray*}

 Since $V(z)$ is not invertible when $\bar{\a} =0,$ it is convenient
to eliminate the last component $p_{\c}(z)$ from the
 {$\c+1$-st} equation provided by the last columns   in
 {\eqref{sysred}:}
 \beq \la{e:last}
    (\Tl_o + \c \Tth \Tm)  {p_\c(z)} = \sum_{i=0}^{\c-1} p_{i}'(z)- \sum_{i=0}^{\c-1}
( \Tl \T  p_a + i \Tth \Tm)p_{i}(z)
%\Eq \Tl_o \pi_{\c,\j}=  \sum_{i=0}^{\c-1} [(\j+1)
%\pi_{i,\j+1}- \Tl \bar p \pi_{i,\j}]
\eeq

The, because only
the last of the first $\c$ equations contains $p_{\c}(z)$, in the
form
{
\begin{eqnarray*}
    (z-\bar p) p_{\c-1}'(z)-p p_{\c-2}'(z) &=& \c \Tmu (\bar \theta + \Tth z) p_{\c}(z)+ \cdots \\
    &=& \fr{\c \mu(\bar  \th +\Tth  z) }{\l_{ob} + \c \Tth \m} \left ( \sum_{i=0}^{\c-1} p_{i}'(z)-
 \sum_{i=0}^{\c-1} ( \Tl \bar p + i \Tth \Tm) p_{i}(z) \right ) + \cdots
\end{eqnarray*}  } we arrive, letting \beq \la{redv}
\widetilde{\vc{p}}(z)=( {p_0(z),} \ldots, p_{\c-1}(z))\eeq denote the
first $\c$ unknowns, to the simplified system:
              \beq \la{redsys} && {\widetilde{\vc{p}}'(z) (V_{\c-1}(z)-
          \k \bff L ) =
          \widetilde{\vc{p}}(z) (U_{\c-1}(z)- \k (\Tl
           \T  p_a \bff L + \Tth \Tm \bff{l}_1)}
\end{eqnarray}
where $V_{\c-1}$  {and} $U_{\c-1}$ are projections of {$V(z)$
and $U(z)$} on the first $\c$ coordinates, where $\bff L $ denotes
a matrix with ones on the last column and $0$ {else,}
$\bff{l}_1$ denotes a matrix with $0,1,2, \ldots$ on the last
column and $0$ else,  {and}
  \beq\la{e:k} \k=\k(z)=\fr{\c
\mu(\bar  \th +\Tth z) }{\l_{ob} + \c \Tth \m}=\fr{\bar  \th +\Tth
z}{\r + \Tth} \eeq

 Explicitly, we have:
 \beq \la{sysred1} &&V(z)=\left[\begin{array}{ccccc}
   z  -\bar p &  -p  &           &    &-\k    \\
       &      z  -\bar p  &   -p   &  & \vdots      \\
& &   \ddots  &   \ddots    &      \vdots \\
   && &    z -\bar p& -p -\k  \\
   &&&&  z -\bar p-\k
\end{array}\right]
\eeq  {\scriptsize {
   \begin{eqnarray*}  U(z)=
\left[\begin{array}{ccccc}
\Theta_0 & \Tl  p_a &   &   & - \k \Tl \T  p_a   \\
\Tmu (\bar  \th +\Tth  z)   & \Theta_1 &  \Tl  p_a  &  &  - \k(\Tl \T  p_a + \Tth \T \mu)\\
  &  2 \Tmu  (\bar  \th +\Tth  z) & \Theta_2 & \Tl  p_a  &  - \k(\Tl \T  p_a + 2 \Tth \T \mu) \\
  &   &  \ddots  & \ddots & \Tl  p_a - \k(\Tl \T  p_a + (\c-2) \Tth \T \mu) \\
  & & & (\c-1) \Tmu (\bar  \th +\Tth  z)  & \Tl( z \T  p_a-1) - (\c-1)\Tmu \bar  \th - (\c-1)\Tmu \Tth (\k+1) -\k \Tl \T  p_a
\end{array}\right ]
\end{eqnarray*}} }

\begin{remark} \la{r:si} The determinant of the reduced system %--see Remark \ref{r:red}
 is
$$Det(V(z))=(z-\bar{p})^{\c-1} (z-\bar{p}-\k(z))=
(z-\bar{p})^{\c-1} \fr{\r}{\r + \Tth} (z-\bar{p}-\Tr (1 +
\fr{\Tth}{\bar \th} \bar p ) )=0$$
  {where  $\Tr=\frac{\c \m \bar
\th}{\l_{ob}}$.}
 { \textbf{Therefore, for $\c
\geq 3$, $\bar{p}$ is an irregular singularity.}} \end{remark}

\begin{remark} %\la{r:st} %Remark 6.2
  We  consider here only the stable  case, when stationary probabilities exist. As is well known, {when
$K=0$, the system is stable  if $ \bar{\a}
>0$ (nonpersistent retrials).} When $ \bar{\a}
=0$, we conjecture that ergodicity holds precisely when the dominant regular singularity
 is outside the unit circle
 \beq \la{e:stab}
    z_r:=\bar{p}+  \fr {1 +  \fr{\Tth}{\bar \th} \bar p }{\r} >1 \Longleftrightarrow    p \r <{1 +  \fr{\Tth}{\bar \th} \bar p }
 \eeq
 where
 \beq \la{e:r} \r:=\fr{\l_{ob} }{\mu\bar \th \c}=\Tr^{-1}.\eeq

 Equivalently,
 \begin{eqnarray*} \x:=p( \r {\bar \th} +  {\Tth})+ \th < 1 \end{eqnarray*}

 This reduces to  Hanschke's condition $p \r < 1 $ when $\Tth=0,$ and to the Atencia-Moreno condition $ \r < 1-\th =\bar \th $ when $\Tth=0$  {and} $p=1$.

 %or by using   Falin's approximation, which % replaces
%the retrial rates by $\I$ when $ \j >L$, changing %the orbit into a regular  queue.

To relate to  \citet{hanschke1999matrix},  note that he assumes $p
+ \bar p \leq  1$  {and} $\a + \bar \a \leq 1,$ while we assume $p
+ \bar p =  1$ (which allows us to use Pearce's Lemma). Our
assumption is w.l.o.g., however to compare with \citet[Thm
1]{hanschke1999matrix}, we must rewrite his ergodicity condition
in our notation as \beq \la{e:Herg} \fr{\l_{ob}}{\c \mu} \n_{q} <
\n_{q} + \n_{a} \eeq which shows that our ergodicity condition
\eqref{e:stab} generalizes his.

Unfortunately, we do not have a reference for the
 rather fundamental and plausible conjecture that
 the obviously necessary ergodicity condition that the dominant singularity lies  outside the unit circle is  also sufficient\fn[4]{We believe this since Kolmogorov's system encodes in principle all the necessary information on our problem,   but we have been unable to check  this beyond the  particular cases  studied by \citet{hanschke1999matrix} and \citet{atencia2003queueing}.}.
 \end{remark}

\begin{remark}  The
  asymptotic behavior may be determined by expanding {the
  function}
   around the  (regular)
   singularity  $z_r$  {---} see Section~\ref{s:asy}.

   Note that the parameters $p$, $\th$ and $\a_0$ play a crucial role in the location of the singularities, unlike $ p_a$.
   \end{remark}

\begin{example} For $K=0$ and $\c=1$ or $\c=2$ servers, the matrices $V$
 {and} $U$ are  {given, respectively, by} $V=z  -\bar p
-\k(z)$  {and} $U=\T \l( z \T p_a-1- \k(z) \T p_a)$, and
\begin{align*}
&&V(z)=\left[\begin{array}{cc}
   z  -\bar p &  -p -\k(z)    \\
   0    &      z  -\bar p -\k(z)  \\
\end{array}\right]\\\no
    &&U(z)=\left[\begin{array}{cc} \Tl( z \T  p_a-1)& \Tl( p_a -\k(z) \T  p_a)\\
\Tm(\bar  \th +\Tth  z)& \Tl( z \T  p_a-1-  \k(z)    \T  p_a)-\Tm
(\bar  \th +\Tth + \k(z)\Tth)             \end{array}\right]
\end{align*}

\end{example}

\begin{example} The solution for $\c=1$ is  {given}, after putting \begin{eqnarray*}
\x:=p( \r +\Tth)+ \th, \; u=1- \x +\rho - \rho z=\bar \th +\bar p
\Tth  - \r(z -\bar p), \bar \l=\Tl (1-\bar p \T  p_a)(1+
\fr{\Tth}{ \r})
    \end{eqnarray*}
 {by}
 \beq \la{e:s=1} &&q(z)= c
   u^{ -\bar \l} e^{-\fr{\Tl \T  p_a}{\r}  u
   }
    \\&&p(z)= \r(1 -\bar p \T  p_a) u^{-1} q(z) \no
\eeq
   with singularity at $z_r >1$ when $\x <1$,  {on the branch cut} $[z_r, \I)$.
   Using  $q(z)+ p(z)|_{z=1}=1$ yields easily the
             proportionality constant $c$.

From \eqref{e:s=1} it may be shown that  asymptotically
             $$(p_{0,\j},p_{1,\j}) =c  (1-\x)^{\bar \l+1} \fr{(\j+1)^{(\bar \l -1)} \x^\j }
             {\Gamma(\bar \l)}
             \left(1,
             \fr{\x (\bar \l+\j)}
             {\bar \l} \right )
\sim_{n \to \I} c (1-\x)^{\bar \l+1}  \fr{\j^{\bar \l-1}}
             {\Gamma(\bar \l)} \x^\j
             (1,\fr{\j }{\Tm })  $$
where $\bar \l=\Tl  (1+\frac{\Tth }{\rho }),$ and $(\j+1)^{(\bar \l -1)}=(\j+1)(\j+2)\cdots$ is an ascending Pochammer symbol, in line with \citet{kim2012tail} (who worked in the classic case with $\x=\r$).

    With $\T  p_a=0,$  the exact solution is:
           \begin{eqnarray*}
           (q(z), p(z))=
        \frac{  1- \x }{  1- \x +\rho} \le(\fr{ 1- \x }{ 1- \x +\rho -\rho z
   }\ri)^{\Tl  (1+\frac{\Tth }{\rho })}
         (1,(\bar \th +\bar p \Tth) \Tr  - (z -\bar p))^{-1})
         \end{eqnarray*}
          a negative binomial distribution.

Furthermore, if  $\bar \th=1$  {and} $p=1,$ this yields
  the classic
             \begin{eqnarray*} (q(z),p(z))=(1-\r) \left(\fr{1-\r}{1-\r z} \right)^{\Tl}
             \left(1, \fr{\r}{1-\r z}\right )
             \end{eqnarray*}
             and $(p_{0,\j},p_{1,\j})= (1-\r)^{\Tl+1}
             \big (\fr{\r^\j \Tl_{(\j)}}
             {\j!},
             \fr{\r^{\j+1} (\Tl+1)_{(\j)}}
             {\j!} \big )$
             see --- \citet{falin1997retrial} pg. 101.
  \end{example}

\section{The classic persistent retrial model ($ \bar \a  =0$) with full acceptance ($p_a=1$) and no feedback to orbit
($\Tth =0 $): a %reduced
generalized Okubo %differential
system \la{s:Oku}}

In this section we consider the  case   {of} $\bar \a=  \Tth=\T
p_a =0$,  {for which} the system \eqref{redsys} {becomes:}
 \beq \la{Ok} &&\widetilde{\vc{p}}'(z) (V_{\c-1}(z)-
          \Tr \bff L ) =
          \widetilde{\vc{p}}(z) U_{\c-1}   %- \Tr \Tl \T  p_a \bff L )}
\end{eqnarray}
where $\Tr =\fr{\c
\mu}{\l_{ob}}$, and
 {$V_{\c-1}(z)$ and $U_{\c-1}$ are given, {respectively, by:}
  \begin{eqnarray*}  &&\left[\begin{array}{cccccc}
   z -\bar p  &  -p  &    &    &   &  -\Tr   \\
       &      z -\bar p   &   -p   &  &  &     -\Tr   \\
& &   \ddots  &   \ddots    &    &    \vdots \\
& & &   z -\bar p & -p &  -\Tr   \\
&&&&  z -\bar p & - p -\Tr \\
&&&&&    z -\bar p -\Tr\\
\end{array}\right]
 \end{eqnarray*}
 \beq \la{e:Ul}
\left[\begin{array}{cccccc}
-\Tl& \Tl  &  \\
 \Bmu  & - \Tl - \Bmu   & \Tl &  \\
 & 2 \Bmu     & -\Tl - 2\Bmu  & \Tl &  \\
& &  \ddots & \ddots & \ddots  &  \\
&& & (\c-2) \Bmu     & -\Tl - (\c-2)\Bmu  & \Tl \\
&& & &(\c-1) \Bmu  &-\Tl - (\c-1) \Bmu  \end{array}\right]
   \eeq
}
where $\Bmu=\Tmu \bar \th$.

   \begin{remark}
For any number  {$\c$ of servers,} $U_{\c-1}=B_{\c-1}- \Tl E,$
where $E$ is a matrix with a $1$ in the lower right corner and $0$
else, is the generator of the time until blocking of an Erlang
loss system with $\c$ servers.
%This turns out to imply, cf. \citet{cohen1957basic},  that the limiting
%marginal distribution of occupied servers under several limiting
%regimes (when $p \to 0$ and  when $\n  \to 0$) is that of an Erlang
%loss system $M(\l)/M(\m )/\c/\c$.
\end{remark}

   Since $\bff U=U_{\c-1}$ is now a constant matrix,
 the reduced system  $ {\widetilde{\vc{p}}'(z)} V_{\c-1}(z) =
  {\widetilde{\vc{p}}(z)} \bff U$ is  of ``generalized Okubo type" defined
by a canonical form\fn[4]{
%\begin{comment}
There is also an equivalent formulation %of  this case
%via the  rank zero "Birkhoff canonical form"
$$\vc{q'}(z)=-\vc{q}(z) \Big( \bff T
+\fr{{\bff U}+ \bff I}{z}\Big)$$  {from which} \eqref{Oku} is
obtained by Laplace transform.
%In the generalized Okubo case,
%higher rank systems may appear by inverse Laplace transform?
%\end{comment}
} : \beq
&& {\widetilde{\vc{p}}'(z)\big(z \bff I- \bff T \big)=\widetilde{\vc{p}}(z) {\bff
U}}  \la{Oku}\eeq with $\bff T$ given by
{  \beq \la{Okmat} &&\left[\begin{array}{cccccc}
   \bar p &  p  &   &  &    &  \Tr  \\
  &      \bar p  &   p   &  &  &     \Tr   \\
&  &   \ddots  &   \ddots    &  &      \vdots \\
&  & & \bar p   &    p&   \Tr  \\
&  & & & \bar p   &   p+  \Tr  \\
 & & &   &  & \bar p+   \Tr \\
\end{array}\right]
\eeq
}

\begin{remark}
The recently introduced class of generalized Okubo-type systems
constitutes a powerful generalization of classic hypergeometric
equations (see for example \citet{balser2006systems,hiroe2009twisted}, and \citet{oshima2011fractional} for an
introduction).  The Jordan form
of $\bff T$ turns out to be useful for classifying
differential systems.

Note  that we have only two Jordan blocks (for any number $\c$ of
servers), for which the eigenvalues $\bar p $ and
$\bar p + \Tr$ of $\bff T$ yield the singular points of our system
(irregular and regular, respectively), and that the Jordan form of
our $\bff T$ is not diagonal.

 {According to \citet{kohno1999global}, Ch4,} a Fuchsian
equation may be written as a system of
 type \eqref{Oku} iff $T$ is diagonal (this case includes several important special functions,
 like the generalized hypergeometric, and
 the {Pochammer system}).% (obtained with
% $U_{i,j}=\l_i - \r(1-\d_{i,j}),$
%when $\sum_{i=1}^n \l_i \neq n \r$).

 Therefore,  the case of persistent retrial queues with $\c \geq 3$ is  qualitatively different from the
   cases $\c=1$ and $\c=2$,  since
     the point $z=\bar p$ is now  a non-Fuchsian irregular singular point,
with  non-diagonal $\bff T,$  a situation that has been much less studied in the literature.
\end{remark}

\begin{remark} When $p >0,$   {putting $y=\fr{z-\bar p}{p}$ reduces} the
problem effectively to the {\bf pure retrial} case {with $p=1$
and} $\bar p=0$, after replacing $\Tr$ and $ \r$ by $\bar
\r=\Tr/p$ and $\x=p \r$, {respectively.} We will call {this
case} the {\bf standardized Okubo problem}. \end{remark}

Recall that we  suppose throughout  that the dominant regular singularity $\bar{p} + \tilde{\r} $
 is outside the unit circle, i.e.   $p <\tilde{\r}$, or simply
$ \tilde{\r} > 1$ in the  standardized case.

In the  {following subsections, we will revisit the cases of
$\c=2$ and $\c=3$, respectively.}

\ssec{Two  servers \label{s:two}}

We note first that {for $\mathbf{\c=2}$}, it is also possible to treat the general Pearce system \eqref{sysred}:
               \begin{eqnarray*} \mathbf{p'}(z)        \left[\begin{array}{ccc}
%\begin{equati\ddotson*}
%\begin{pmatrix}
   z -\bar p  &  -p  &     1         \\
   0    &      z  -\bar p  &   1      \\
   0    &   0    &     {\bar{\a}  } \\
\end{array}\right]=\mathbf{p}(z) \left[\begin{array}{ccc}
%\begin{equati\ddotson*}
%\begin{pmatrix}
\Tl (\T  p_a z -1) &     \Tl p_a  &          \Tl \T  p_a \\
       \Tmu (\bar \theta+ \T{\th}  z)  \  & \Tl (\T  p_a z -1)- \Tmu \bar \theta  &
       \Tl \T  p_a   +\Tmu  \T{\th}   \\
       0              &   2 \Tmu (\bar \theta+ \T{\th}  z) & \Tl_{ob} + 2 \Tmu  \T{\th} \\
\end{array}\right]   \end{eqnarray*}

 However, we consider only  the Okubo case (with $\bar \a=0$,  with full acceptance $\T p_a=0$ and no feedback $   \Tth =0$), when
 the  system is:
            \begin{eqnarray*} \mathbf{p'}(z)        \left[\begin{array}{ccc}%\begin{pmatrix}
   z -\bar p  &  -p  &     1         \\
   0    &      z  -\bar p  &   1      \\
   0    &   0    &   0   \\\end{array}\right]=\mathbf{p}(z) \left[\begin{array}{ccc}
%\begin{equati\ddotson*}
%\begin{pmatrix}
-\Tl  &     \Tl               &          0               \\
       \Bmu   \      & -\Tl - \Bmu   &    0              \\
       0              &   2 \Bmu  & \Tl_{ob}  \\
\end{array}\right]   \end{eqnarray*}
(compare with \citet[(4.3),(4.4)]{hanschke1987explicit}, who assumes
 $\bar p=0$, $\bar \th=1$,  {and}  $\l_{ob}=\l$).

After eliminating $p_2(z)$ using \eqref{e:last}, the reduced
Okubo system  \eqref{Oku}  {and} \eqref{Okmat} for  $\widetilde{\vc{p}}(z)=(p_0(z),p_1(z))$ becomes:\fn[4]{\citet{hanschke1987explicit}
obtains when $p=1$ classic hypergeometric solutions. He considers
the second order scalar equation $L q(z)=0$ for $q(z)=p_0(z),$
which is:
             \begin{eqnarray*} &&L= (\Tl+ \Bmu +(z-\Tr) D) (\Tl +
           z D)- \Bmu (\Tl+  (1+\Tr) D  )=\\&&\Tl ^2   +   (z   (2 \Tl + \
\Bmu + 1 ) - \Tr  (3 \Tl/2 +\Bmu + 1 )) D +
    z (z  - \Tr )  D^2\end{eqnarray*}
           (using $D z D=D + z D^2$), where $\Tr=\fr{\mu \c}{\l_{ob} }$ was defined in \eqref{e:r}.
Hanschke  \citeyearpar[(4.8)]{hanschke1987explicit}  notes  that
    putting $\r z=x$ shifts the singularity $z=\r^{-1}$, yielding
    the canonic Gauss hypergeometric equation
    $$x(x-1) q''(x)+[x(2\Tl +\Bmu +1)- (\frac{3 \Tl}{2}+\Bmu +1)]
    q'(x)+ \Tl^2 q(x)=0$$
    whose only analytic solution in the unit disk is the
    Gauss hypergeometric series. This determines all unknowns
    up to a  proportionality constant, obtained using
    $q(z)+ p_1(z)+ p_2(z)|_{z=1}=1$.
The explicit formulas for $p_i(z)$ are given in \citet[Thm.
2.1]{falin1997retrial},  {which lead to}  some   explicit
formulas
    for classic performance measures,  {such as $E[N]=\sum_{i=0}^2 p_i'(1)$, given in} \citet[Thm. 2.2]{falin1997retrial}
    (see also \citet[(91-93)]{falin1990survey}). Further factorial moments
    $\sum_{i=0}^2 p_i''(1)=E[N(N-1)]$, $\sum_{i=0}^2  p_i'''(1)=E[N(N-1)(N-2)]$,
    \ldots
    are given in the last section of \citet{hanschke1987explicit}.}
 \beq \la{e:O2} &&\widetilde{\vc{p}}'(z) \left[\begin{array}{cc}
   z  -\bar p &  -p -\r    \\
   0    &      z  -\bar p -\r  \\
\end{array}\right]= \widetilde{\vc{p}}(z) \left[\begin{array}{cc} -\Tl& \Tl\\
\Tm& -\Tl-\Tm         \end{array}\right]\eeq

\sec{Shifting to the dominant regular singularity and asymptotics for the pure retrial case
%of the stationary probabilities
 \la{s:asy}}
%\ssec{The Birkhoff canonical form \la{s:Birk}}

  Consider now  the standardized Okubo problem
  $\widetilde{\vc{p}}'(y)\big(y \bff I- \bff T \big)=\widetilde{\vc{p}}(y) {\bff
U}$ {
  with $U$ defined in \eqref{e:Ul}, and $\bff T= T_+ +\bar \r \bff L$,
   {where}
\[
  T_+=\left(
\begin{array}{cccccc}
 0 & 1 &  \\
 & 0 & 1&  \\
  & & 0& 1 & \\
   & & &  \ddots&\ddots& \\
   &&&& \ddots
\end{array}\right), \q \bff L=\left(
\begin{array}{cccccc}
 0 &0 & \dots &1 \\
  0 &0 & \dots &1 \\
   \vdots &  \vdots &  \vdots&\vdots \\
\end{array}\right)
\]
  in which the singularities of the system are the eigenvalues  $0$, which is irregular, and $\bar \r$, which is regular.}

   Since $\bff T^{\c}=\bar \r
\bff T^{\c-1},$ the resolvent $\big(  y \bff I- \bff T \big)^{-1}$ may be
easily computed, yielding:
{
\begin{eqnarray*} &&( y I-\bff T)^{-1}=y^{-1} \left(\bff I + \fr{1}{y } \bff T
 + \cdots + \fr{1}{y^{\c -2} }
\bff T^{\c -2}+ \fr{1}{y^{\c -2}(y-\bar \r )} \bff T^{\c -1}\right)
\end{eqnarray*}

For $\c=3$ for example, we find:
{
\begin{eqnarray*} ( y I-\bff T)^{-1} &=& y^{-1} \left(\bff I + \fr{1}{y}
\bff T + \fr{1}{y(y-\bar \r )} \bff T^2 \right) =y^{-1} \left[\bff I + \fr{1}{y }
\bff T + \fr{1}{\bar \r } \left(\fr{1}{y-\bar \r}-\fr{1}{y}\right)\bff T^2\right] \\
&=& y^{-1} \left[\bff I + \fr{1}{y }\left(\bff T -\fr{1}{\bar \r } \bff T^2\right)+ \fr{1}{\bar \r } \fr{1}{y-\bar \r}\bff T^2 \right] \\
&=& y^{-1} \left[\bff I + \fr{1}{y } \bep 0&1&-(1+\x) \\0&0&0\\0&0&0 \eep+
 \fr{1}{y-\bar \r } \bep 0&0&\x  + 1+ \x ^{-1} \\0&0&1+ \x ^{-1}\\0&0&\x ^{-1} \eep\right]
\end{eqnarray*}   }
where we used $\bff L  ^2=\bff L $, $\bff L  T_+=\bff 0$, $T_+ \bff L  :=\bff L _1=\left(
\begin{array}{lll}
 0 & 0 & 1 \\
 0 & 0 & 1 \\
 0 & 0 & 0
\end{array}\right),$ and $\bff T^2=\bep 0&0&\x  + 1+ \x ^{-1} \\0&0&1+ \x ^{-1}\\0&0&\x ^{-1} \eep$.

Multiplying} the resolvent on the RHS, we get  $\forall \c \geq 3$ a  differential system
$$ \widetilde{\vc{p}}'(y)= \widetilde{\vc{p}}(y) \bff R=\widetilde{\vc{p}}(y) \le((y-\bar \r)^{-1} \bff D + y^{-1} \bff D_1 + y^{-2} \bff D_2\ri)$$ of Poincar\'e  rank $1$ %$(\c-2)_+$
at $y=0$ (thus this   singularity is nonregular) and of rank $0$ at $y=\bar \r$,
with \beq && \bff D= \left(
\begin{array}{ccc}
  & 0 & \vdots \\
  & &  0  \\
   & 0&    \l(1+\r^{-1})\\
   &0& -\l -(\c -1) \mu
\end{array}\right),   \\&&\bff D_1= \left(
\begin{array}{cccc}
  -\l& 0 & \dots \\
  0& -\l -\mu &  0  \\
   & &  \ddots  &\\
&& -\l -(\c -2) \mu&-\l(1+\r^{-1})\\
   0&0&0&0
\end{array}\right), \bff D_2= \left(
\begin{array}{cccc}
 0  &\l &0 &\\
  0& 0  &\l&0\\
   & & 0&  0\\
   && &
\end{array}\right)  . \la{firdec} \no%\\\Eq&& y^{\c-2}  \widetilde{\vc{p}}'(y) = \widetilde{\vc{p}}(y) \left(D_1 +
 % \T{B}_0 y^{-1} +     \fr{1}{\bar \r- y }D\right) %\la{secdec}
  \eeq

The asymptotics of the stationary probabilities, and in fact a full asymptotic expansion,  may be obtained  by looking for logarithmic free
regular solutions  $\widetilde{\vc{p}}(y)= (\bar \r-y)^w \vc{g}(y),$ where $\vc{g}(y)$ is analytic at the regular singularity
$y=\bar \r.$  The eigenvalue $w$, as well as the first term of such an expansion, are obtained from the spectral decomposition of the matrix $\bff D$.

In principle, this expansion requires an application of the  vector version of Frobenius's  formal series
 expansions around singular points --- see for example \citet[Ch 4]{coddington1955theory}, \citet[Ch 4]{wasow2002asymptotic}, and, for recent developments,  \citet{barkatou1999algorithm}.

 Indeed, the first term of such an expansion was recently provided by \citet{liutail} and \citet{kim2012tail}
 (who use a different dimensionality reduction from $\c+1$ to $\c$, by the important probabilistic technique of  {``}censoring").

  Note however that the    vector  Frobenius method
 is available in several flavors (which typically require
$\bff D$ to be non-singular, which is not the case in  our
differential system, nor in  {\citet{liutail} or} \citet[Thm
1]{kim2012tail}), and therefore obtaining a full asymptotic
expansion requires further work.

{\bf Acknowledgements:} We thank Jesus Artalejo, JL.
Lopez,  A. Quadrat,  T. Cluzeau, F. Beukers,  C. Koutschan and the anonymous referees for very useful advice. %\end{comment}

The first two authors were partially supported by the
French-Romanian Programme LEA Math-Mode.

\textbf{F. Avram} is a full professor in the Mathematics department at  Universit\'e de Pau, France. He received his Ph.D. degree in Mathematics
from Cornell University in 1986. He has published over 60 research papers in  internationally-refereed journals, in the field of applied probability and stochastic processes (with particular emphasis on queueing networks, risk theory and mathematical finance).

\textbf{D. Matei} is a senior researcher in the Geometry and Topology group at the ``Simion Stoilow''
Institute of Mathematics of the Romanian Academy in Bucharest. He received his Ph.D. degree in Mathematics
from Northeastern University, Boston in 1999. He has published research papers in algebraic geometry and topology.
His research interests include hypergeometric differential equations,
birth-and-death processes and random walks.

\textbf{ Y.Q. Zhao} received his Ph.D. degree at the  University of Saskatchewan. He has been a full professor at Carleton University, Ottawa, since 2003, and was the Director of the School of Mathematics and Statistics from 2004 to 2007 and from 2010 to 2011.
Dr. Zhao's research interests are in applied probability and stochastic processes, with particular emphasis on computer and telecommunication network applications. He has published over 70 papers in internationally-refereed journals.

\bibliographystyle{plainnat}

\bibliography{ret}

\end{document}